\documentclass{file}
\usepackage{pat}
\usepackage{paralist}
\usepackage{dsfont}  
\usepackage[utf8x]{inputenc}
\usepackage{skull}
\usepackage{paralist}
\usepackage{centernot}
\usepackage{mathtools}
\usepackage{thmtools}
\usepackage{hyperref}
\hypersetup{hidelinks}
\usepackage{needspace}
\usepackage[table]{xcolor}
\usepackage{stmaryrd,tikz}
\usetikzlibrary{calc}
\usepackage{graphicx}
\usepackage{booktabs}
\usepackage{listings}
\usepackage{pifont}%
\usepackage[normalem]{ulem}
\usepackage{cancel,float}
\usepackage{comment}

\usepackage{todonotes}
\usepackage{mathdots}
\usepackage{bm}
\usepackage{etoolbox}
\definecolor{maroon}{rgb}{0.5, 0.0, 0.0}
\definecolor{darkblue}{rgb}{0.0, 0.0, 0.55}

\usepackage[longnamesfirst,numbers,sort&compress]{natbib}

\usepackage[mathlines]{lineno}
\newcommand*\patchAmsMathEnvironmentForLineno[1]{%
 \expandafter\let\csname old#1\expandafter\endcsname\csname #1\endcsname
 \expandafter\let\csname oldend#1\expandafter\endcsname\csname end#1\endcsname
 \renewenvironment{#1}%
    {\linenomath\csname old#1\endcsname}%
    {\csname oldend#1\endcsname\endlinenomath}}%
\newcommand*\patchBothAmsMathEnvironmentsForLineno[1]{%
 \patchAmsMathEnvironmentForLineno{#1}%
 \patchAmsMathEnvironmentForLineno{#1*}}%
\AtBeginDocument{%
\patchBothAmsMathEnvironmentsForLineno{equation}%
\patchBothAmsMathEnvironmentsForLineno{align}%
\patchBothAmsMathEnvironmentsForLineno{flalign}%
\patchBothAmsMathEnvironmentsForLineno{alignat}%
\patchBothAmsMathEnvironmentsForLineno{gather}%
\patchBothAmsMathEnvironmentsForLineno{multline}%
}

\newcommand{\F}{\mathbb{F}}
\newcommand{\rank}{\operatorname{rank}}
\newcommand{\tr}{\operatorname{tr}}

\newcommand{\diag}{\operatorname{diag}}

\newcommand{\one}{\mathbf{1}}
\newcommand{\eps}{\varepsilon}

\newcommand{\ind}{\mathbf{1}}

\definecolor{brightmaroon}{rgb}{0.76, 0.13, 0.28}
\definecolor{linkblue}{rgb}{0, 0.337, 0.227}
\newrobustcmd{\onesub}{\mathord{\includegraphics{one-sub}}}
\newrobustcmd{\leftup}{\mathord{\includegraphics{left-up}}}
\makeatletter
\newcommand{\xMapsto}[2][]{\ext@arrow 0599{\Mapstofill@}{#1}{#2}}
\def\Mapstofill@{\arrowfill@{\Mapstochar\Relbar}\Relbar\Rightarrow}
\makeatother

\title{\MakeUppercase{Dual near-factorizations and Schur products in association schemes}}

\author{
Bobby Miraftab
\thanks{School of Computer Science,
Carleton University.}
\quad
Sho Suda\thanks{Department of Mathematics, National Defense Academy of Japan, Yokosuka, Kanagawa 239-8686,
Japan.}}

\date{}

\begin{document}

\maketitle

\begin{abstract}
Let $G$ and $H$ be Gram matrices of spherical $2$-designs on the same number of indices. 
We study when their Schur product is again the Gram
matrix of a spherical design. 
We give a criterion in terms of orthogonal projections with constant diagonal, and an explicit $15$-point spherical $3$-design whose Gram matrix belongs to no symmetric Bose--Mesner algebra.
For a symmetric association scheme, the factorization problem takes the form $v(E_S\circ E_T)=\mu E_U$. This is a counterpart of ordinary matrix
product factorizations in the adjacency basis. The case $E_U=I-E_0$ gives a dual near-factorization and has a regular simplex as its normalized
output. We characterize complementary input supports by a strongly regular fusion, give a reciprocal criterion, and classify the case in which one factor is the last Hamming idempotent. 
\end{abstract}

\noindent\textbf{Keywords:}
Association schemes, Krein parameters, Schur products,
spherical designs, dual near-factorizations,
primitive idempotents, strongly regular graphs.

\medskip
\noindent\textbf{2020 Mathematics Subject Classification:}
Primary 05E30; Secondary 05B30, 05C50, 52C35.

\section*{Generative AI disclosure}

OpenAI's ChatGPT and Anthropic's Claude were used during the development of this work to assist with language editing, literature searches, exploration of examples, and the formulation and checking of some computational arguments. AI-generated suggestions were not accepted without verification. All results, proofs, computations, bibliographic information, and claims in the final manuscript were checked by the authors, who take full responsibility for the content.

\section{Introduction}
Let $G$ and $H$ be the Gram matrices of two spherical $2$-designs on
the same number of indices. Their Schur product $G\circ H$ is the Gram
matrix of the pointwise tensor products of the two configurations. We
ask when these tensors form another spherical design in their linear
span. Gram matrices and tight frames provide a way to study this
question; see
\cite{delsarte-goethals-seidel,benedetto-fickus,peng-waldron,waldron-gramian}.
Throughout the paper, configurations are indexed collections of unit
vectors, and repetitions are allowed unless distinctness is stated.

We first express the spherical $2$-design condition in terms of
orthogonal projections with constant diagonal. The framework does not
require an association scheme. Example~\ref{ex:twisted-pentagons} gives
two planar input configurations whose Gram product represents fifteen
distinct nonantipodal points forming a spherical $3$-design in $\R^4$.
Its output Gram matrix belongs to no symmetric Bose--Mesner algebra.
We also give a third-moment criterion and prove that pointwise tensor
products with both input dimensions at least two cannot have strength
four, even when their span is a proper subspace of the tensor product.

Herman and Miraftab \cite{herman2026association} study ordinary matrix product factorizations
\[
 A_SA_T=A_U,\qquad A_S=\sum_{i\in S}A_i,
\]
where the factors and output are loopless unions of relations in an association scheme. 
More generally, matrix product factorizations $A(G)=A(H)A(K)$ have been studied from structural
\cite{maghsoudi2025matrixproduct}, spectral \cite{akbari2025spectral}, and primeness \cite{akbari2025prime} perspectives, as well as for infinite graphs \cite{miraftab2026infinite},
Cayley graphs \cite{herman2025cayley}, and
approximate factorizations \cite{maghsoudi2026epsilon}.

The Schur-product question has a corresponding formulation in the primitive-idempotent basis. 
If $E_0,\ldots,E_d$ are the primitive idempotents of a symmetric scheme on $v$ vertices, put $E_S=\sum_{i\in S}E_i$. 
We study identities $v(E_S\circ E_T)=\mu E_U$, where  $\mu>0$, which we call \defin{Krein factorizations}. 
After scaling the basis to $F_i=vE_i$, the intersection numbers in ordinary multiplication are replaced by the Krein parameters in Schur multiplication. The case $\mu=1$ is the direct counterpart of the matrix product factorization
above. 
Allowing a scalar is natural for the Gram problem: if the three
idempotent sums have ranks $r,s,t$, then $rs=\mu t$.

Near-factorizations in association schemes, of the form $A_SA_T=\lambda(J-I)$, are studied by Herman, Lacaze-Masmonteil, Meagher, and Monterde \cite{herman2026near}. 
Their counterpart is
\[
 v(E_S\circ E_T)=\mu(I-E_0),\qquad E_0=J/v.
\]
We call this a \defin{dual near-factorization}. 
The term refers to the
two multiplication structures in the Bose--Mesner algebra; a combinatorial dual is not assumed to exist. For nonprincipal factors, the normalized output is necessarily a regular simplex. The questions in this case concern its factors and the restrictions on the ambient scheme.

\section{Gram products of spherical designs}\label{sec:projections}

\subsection{Gram matrices and constant-diagonal projections}

A finite indexed collection of unit vectors in $\R^r$ is a spherical
$t$-design if the average of every polynomial of degree at most $t$
agrees with its average over the unit sphere. 
We allow repetitions in
an indexed collection. In particular, $x_1,\ldots,x_v\in\R^r$ form a
spherical $2$-design if and only if
\[
 \sum_{i=1}^v x_i=0,\qquad
 \sum_{i=1}^v x_ix_i^{\mathsf T}=\frac vr I_r.
\]
We always take the dimension of a configuration to be the dimension of
its linear span.

\begin{thm}\label{thm:gram-2design}
Let $G,H$ be the Gram matrices of two spherical $2$-designs with $v$
points, in dimensions $r$ and $s$, respectively. Put
$K=G\circ H$ and $t=\rank (K)$. Then $K$ is the Gram matrix of a
spherical $2$-design in $\R^t$ if and only if
\begin{equation}\label{eq:gram-2design}
 GH=0,\qquad K^2=\frac vt K.
\end{equation}
Equivalently, we have
\begin{equation}\label{eq:gram-2design-moments}
 \sum_{i,j=1}^vG_{ij}H_{ij}=0,\qquad
 \sum_{i,j=1}^vG_{ij}^2H_{ij}^2=\frac{v^2}{t}.
\end{equation}
\end{thm}

\begin{proof}
The matrix $K$ is the Gram matrix of the unit vectors
$x_i\otimes y_i$, where $G,H$ are the Gram matrices of
$\{x_i\}$ and $\{y_i\}$. A Gram matrix of $v$ unit vectors and
rank $t$ represents a spherical $2$-design if and only if
$K\one=0$ and $K^2=(v/t)K$.
Now
\[
 \one^{\mathsf T}K\one=\sum_{i,j}G_{ij}H_{ij}=\tr(GH).
\]
Since $K$ is positive semidefinite, the left side is zero if and only
if $K\one=0$. Since $G,H$ are positive semidefinite, the right side
is zero if and only if $GH=0$. This proves
\eqref{eq:gram-2design}. Finally, $\tr K=v$, so Cauchy--Schwarz
applied to its $t$ positive eigenvalues gives
$\tr(K^2)\ge v^2/t$, with equality exactly when
$K^2=(v/t)K$. This proves \eqref{eq:gram-2design-moments}.
\end{proof}

We next restrict the rank for the corresponding projections.
We use the Schur product theorem in the form that
$M,N\succeq0$ implies $M\circ N\succeq0$, see \cite{horn-johnson}.

\begin{prop}\label{thm:constant-diagonal}
Let $P,Q\in M_v(\R)$ be nonzero orthogonal projections of ranks $r,s$,
with constant diagonal entries $r/v,s/v$, respectively. Suppose
\begin{equation}\label{eq:PQ-alphaR}
 P\circ Q=\alpha R
\end{equation}
for some $\alpha>0$ and an orthogonal projection $R$ of rank $t$.
Then $R$ has constant diagonal $t/v$, and
\begin{equation}\label{eq:alpha-rst}
 \alpha=\frac{rs}{vt},\qquad
 \max\{r,s\}\le t\le\min\{rs,v\}.
\end{equation}
In particular, $\alpha\ge1/v$, with equality if and only if $t=rs$.
\end{prop}

\begin{proof}
The diagonal entries of $P\circ Q$ equal $rs/v^2$. Taking traces in
\eqref{eq:PQ-alphaR} gives $\alpha t=rs/v$, and hence the stated
diagonal of $R$.
The Hadamard-rank inequality gives $t\le rs$.
Since $0\preceq Q\preceq I$, the Schur product theorem gives
\[
 0\preceq P\circ Q\preceq P\circ I=\frac rv I.
\]
Thus $\alpha\le r/v$, and similarly $\alpha\le s/v$.
Substituting the value of $\alpha$ yields $t\ge s,r$.
The remaining assertions follow from $t\le v$ and the formula for $\alpha$.
\end{proof}

For clarity, the unscaled identity $P\circ Q=R\ne0$ forces
$P=Q=R=I$. Indeed, $\alpha=1$ in the proof forces $r=s=v$.
The scale $1/v$ has a different meaning: it is the case in which the
pointwise tensor products span the full tensor space. To see this,
write $P=VV^{\mathsf T}$ and $Q=WW^{\mathsf T}$ with orthonormal
columns, and form the matrix $Z$ whose $x$th row is $v_x\otimes w_x$. 
Then we obtain
\begin{equation}\label{eq:ZZT}
 P\circ Q=ZZ^{\mathsf T}.
\end{equation}
Under \eqref{eq:PQ-alphaR}, the condition $\alpha=1/v$ is equivalent
to $Z^{\mathsf T}Z=I_{rs}/v$.

\subsection{The third and fourth moments}

\begin{prop}\label{prop:gram-3design}
Let $K$ be the Gram matrix of a spherical $2$-design. It represents a
spherical $3$-design if and only if
\begin{equation}\label{eq:gram-3design}
 \sum_{x,y}K_{xy}^3=0.
\end{equation}
\end{prop}

\begin{proof}
Let $K$ be the Gram matrix of $z_1,\ldots,z_v$. The additional condition
for strength three is $\sum_xz_x^{\otimes3}=0$. The squared norm of this
tensor is $\sum_{x,y}\langle z_x,z_y\rangle^3$, proving the assertion.
\end{proof}

The next obstruction also concerns the actual span of the output.
\begin{thm}\label{prop:no-tensor-4design}
Let $x_1,\ldots,x_v\in\R^r$ and $y_1,\ldots,y_v\in\R^s$ be unit
vectors spanning their respective spaces. If $r,s\ge2$, then the
pointwise tensor products $x_i\otimes y_i$ do not form a spherical
$4$-design in their linear span.
\end{thm}

\begin{proof}
Identify $\R^r\otimes\R^s$ with the space of $r\times s$ matrices,
and let $L$ be the span of the rank-one matrices $x_i y_i^{\mathsf T}$.
Every $2\times2$ minor vanishes at these matrices. 
If the configuration were a spherical $4$-design in $L$, the average of the square of each minor would agree with its integral over the unit sphere of $L$.
The integral would be zero. 
The integrand is continuous and nonnegative, so each minor would vanish on that sphere and, by homogeneity, identically on $L$.

Next we choose $i,j$ such that $x_i,x_j$ are linearly independent.
Since $x_i y_i^{\mathsf T}+x_j y_j^{\mathsf T}$ belongs to $L$
and has rank at most one, $y_i,y_j$ are proportional.
As the $y_i$ span a space of dimension at least two, choose $k$ such that $y_k$ is not proportional to $y_i$.
Applying the same argument to $(i,k)$ and $(j,k)$ forces $x_k$
to be proportional to both $x_i$ and $x_j$, a contradiction.
\end{proof}

\subsection{An example outside association schemes}

\begin{exa}\label{ex:twisted-pentagons}
Let $u(\theta)=(\cos\theta,\sin\theta)$ and
$(\phi_1,\phi_2,\phi_3)=(0,0,\pi/10)$. For $a=1,2,3$ and
$b=0,\ldots,4$, put
\[
 x_a=u\!\left(\frac{2\pi(a-1)}3\right),\qquad
 y_{ab}=u\!\left(\frac{2\pi b}5+\phi_a\right),\qquad
 z_{ab}=x_a\otimes y_{ab}.
\]
The indexed collections $(x_a)_{a,b}$ and $(y_{ab})_{a,b}$ are spherical
$2$-designs with Gram matrices $G,H$. Their Schur product $K=G\circ H$
is the Gram matrix of the $z_{ab}$.
The triangle and pentagon moments give
\[
 \sum_{a,b}z_{ab}=0,\qquad
 \sum_{a,b}z_{ab}z_{ab}^{\mathsf T}=\frac{15}4I_4,\qquad
 \sum_{a,b}z_{ab}^{\otimes3}=0.
\]
Thus the output is a spherical $3$-design in $\R^4$, of strength exactly
three by ~\Cref{prop:no-tensor-4design}. Each fibre is a pentagon,
and inner products between fibres have absolute value at most $1/2$;
hence the fifteen points are distinct and nonantipodal.

Finally, direct calculation gives
\[
 \sum_{c,d}K_{ab,cd}^{\,5}
 =\begin{cases}
   155/512,&a=1,2,\\
   160/512,&a=3.
  \end{cases}
\]
Every Schur power of a matrix in a Bose--Mesner algebra has constant
row sums. Therefore $K$ belongs to no symmetric Bose--Mesner algebra.
\end{exa}

\section{Krein factorizations in association schemes}\label{sec:schemes}

\subsection{The two multiplication structures}

Let $\mathfrak X=(X,\{R_0,R_1,\ldots,R_d\})$ be a symmetric
association scheme on $v=|X|$ vertices. Let $A_0,\ldots,A_d$ be its
adjacency matrices and $E_0,\ldots,E_d$ its primitive idempotents,
ordered so that
\[
 A_0=I,\qquad E_0=\frac1vJ,\qquad
 \sum_{i=0}^dA_i=J,\qquad \sum_{i=0}^dE_i=I.
\]
Each $E_i$ is a real orthogonal projection of rank $m_i$, with
constant diagonal $m_i/v$. The intersection numbers and Krein
parameters are defined by
\begin{equation}\label{eq:krein}
 A_iA_j=\sum_{k=0}^dp_{ij}^kA_k,\qquad
 E_i\circ E_j=\frac1v\sum_{k=0}^dq_{ij}^kE_k.
\end{equation}
The Krein parameters are nonnegative and satisfy
$q_{ij}^0=\delta_{ij}m_i$. We refer to
\cite{bannai-ito,BCN,delsarte} for standard background.
The first and second eigenmatrices have conventions
\begin{equation}\label{eq:eigenmatrices}
 A_i=\sum_{j=0}^dP_{ji}E_j,\qquad
 E_i=\frac1v\sum_{h=0}^dQ_{hi}A_h.
\end{equation}
Thus $P_{0i}=k_i$ is the valency of $R_i$, and $Q_{0i}=m_i$.
The scheme is formally self-dual if compatible orderings make $P=Q$.
Put $F_i=vE_i$. The second identity in \eqref{eq:krein} becomes
\[
 F_i\circ F_j=\sum_{k=0}^dq_{ij}^kF_k.
\]
Thus $A_SA_T=A_U$, studied in \cite{herman2026association}, has the
coefficient-one counterpart $F_S\circ F_T=F_U$. We allow a positive
scalar on the right, as required by the normalization of Gram products.

For $S\subseteq\{0,1,\ldots,d\}$, put
\[
 \begin{aligned}
 A_S&=\sum_{i\in S}A_i,& E_S&=\sum_{i\in S}E_i,\\
 m(S)&=\sum_{i\in S}m_i,& \eta_h(S)&=\sum_{i\in S}Q_{hi}.
 \end{aligned}
\]
For two index sets $S,T$, define the aggregate Krein coefficients by
\begin{equation}\label{eq:aggregate}
 c_k(S,T)=\sum_{i\in S}\sum_{j\in T}q_{ij}^k.
\end{equation}
Then we have
\begin{equation}\label{eq:aggregate-expansion-new}
 v(E_S\circ E_T)=\sum_{k=0}^dc_k(S,T)E_k.
\end{equation}
We call an identity
\begin{equation}\label{eq:projection-valued}
 v(E_S\circ E_T)=\mu E_U,\qquad \mu>0,
\end{equation}
with nonempty supports $S,T,U$, a \defin{Krein factorization}.
Equivalently, $E_S\circ E_T$ is a positive scalar multiple of an
orthogonal projection; we also call such a product \defin{projection-valued}.
In comparison with
\eqref{eq:PQ-alphaR}, the scalar used throughout the scheme setting is
$\mu=v\alpha$.

\subsection{Coefficient and design criteria}

\begin{prop}\label{thm:coefficient-criteria}
Let $S,T,U$ be nonempty subsets of $\{0,1,\ldots,d\}$, and let
$\mu>0$. The following are equivalent:
\begin{enumerate}
\renewcommand{\labelenumi}{(\roman{enumi})}
\item $v(E_S\circ E_T)=\mu E_U$;
\item $c_k(S,T)=\mu\ind_{k\in U}$ for $0\le k\le d$;
\item $\eta_h(S)\eta_h(T)=\mu\eta_h(U)$ for $0\le h\le d$.
\end{enumerate}
Consequently, $E_S\circ E_T$ is projection-valued if and only if all
its positive aggregate Krein coefficients are equal. In that case
$U=\{k:c_k(S,T)>0\}$, and $U$ and $\mu$ are uniquely determined by
the two factors.
\end{prop}

\begin{proof}
The equivalence of (i) and (ii) follows from
\eqref{eq:aggregate-expansion-new} and linear independence of the
primitive idempotents. By \eqref{eq:eigenmatrices} and
$A_h\circ A_j=\delta_{hj}A_h$,
\[
 v(E_S\circ E_T)
 =\frac1v\sum_{h=0}^d\eta_h(S)\eta_h(T)A_h.
\]
Comparison with $\mu E_U$ in the adjacency basis gives (iii).
The last assertion follows from (ii) and nonnegativity of the
Krein parameters.
\end{proof}

\begin{cor}\label{cor:scheme-spherical-designs}\label{cor:multiplicity-restrictions}
Let $S,T\subseteq\{1,\ldots,d\}$ be nonempty, and put
$r=m(S)$ and $s=m(T)$. The entrywise product of the normalized
Gram matrices $(v/r)E_S$ and $(v/s)E_T$ is the Gram matrix of a
spherical $2$-design if and only if $S\cap T=\varnothing$ and
all positive coefficients $c_k(S,T)$ are equal.
If their common value is $\mu$ and
$U=\{k:c_k(S,T)>0\}$, then $0\notin U$ and, with $t=m(U)$,
\begin{equation}\label{eq:normalized-gram-product}
 \left(\frac vrE_S\right)\circ\left(\frac vsE_T\right)
 =\frac vtE_U,\qquad rs=\mu t.
\end{equation}
Moreover,
\begin{equation}\label{eq:scheme-rank-window}
 \max\{r,s\}\le t\le\min\{rs,v-1\},\qquad r+s\le v-1.
\end{equation}
\end{cor}

\begin{proof}
Since $E_SE_T=0$ is equivalent to $S\cap T=\varnothing$, the
criterion follows from ~\Cref{thm:gram-2design} and
~\Cref{thm:coefficient-criteria}. Also,
\begin{equation}\label{eq:principal-general}
 c_0(S,T)=\sum_{i\in S\cap T}m_i.
\end{equation}
Comparison of diagonal entries in \eqref{eq:projection-valued}
gives $rs=\mu t$, and hence \eqref{eq:normalized-gram-product}.
The rank bounds follow from ~\Cref{thm:constant-diagonal}.
Finally, disjointness gives $r+s\le\sum_{i=1}^dm_i=v-1$.
\end{proof}

The coefficient criterion follows by expansion in the two standard
bases. Formal self-duality gives an exact transfer between the two
factorization problems.

\begin{prop}\label{prop:formal-transfer}
Suppose that compatible labelings give $P=Q$. For nonempty supports
$S,T,U$ and $\mu>0$,
\[
 v(E_S\circ E_T)=\mu E_U
 \quad\Longleftrightarrow\quad
 A_SA_T=\mu A_U.
\]
In either case $\mu$ is a positive integer.
\end{prop}

\begin{proof}
In these labelings $p_{ij}^k=q_{ij}^k$. The identities therefore have
the same aggregate structure constants. Since $A_U$ is a nonzero
$01$-matrix and $A_SA_T$ has integer entries, $\mu$ is an integer.
\end{proof}

\begin{cor}\label{prop:scheme-3design}\label{sec:proper-support}
The normalized embedding from a nonprincipal sum $E_U$ has strength
at least three if and only if
\begin{equation}\label{eq:scheme-3design-criterion}
 c_k(U,U)=0\qquad(k\in U).
\end{equation}
When this condition holds and the embedding is obtained from
a Krein factorization with both input ranks at least two,
its strength is exactly three.
\end{cor}

\begin{proof}
The embedding is a spherical $2$-design, and
\[
 \sum_{x,y}(E_U)_{xy}^3
 =\tr\bigl((E_U\circ E_U)E_U\bigr)
 =\frac1v\sum_{k\in U}m_kc_k(U,U).
\]
All terms on the right are nonnegative. Apply
~\Cref{prop:gram-3design} and
~\Cref{prop:no-tensor-4design}.
\end{proof}

\subsection{Some structural cases}

Multiplicity-one factors give a simple source of factorizations.
We include the argument to distinguish these from the examples with
both input ranks greater than one.

\begin{prop}\label{thm:sign-conjugation}
If $m_i=1$, then there is an involution $\sigma_i$ of
$\{0,1,\ldots,d\}$ preserving multiplicities such that
\begin{equation}\label{eq:permutation-idempotents}
 v(E_i\circ E_j)=E_{\sigma_i(j)}\qquad(0\le j\le d).
\end{equation}
It satisfies $\sigma_i(0)=i$ and $\sigma_i(i)=0$.
Consequently, $v(E_i\circ E_T)=E_{\sigma_i(T)}$ for every support $T$.
\end{prop}

\begin{proof}
Since $E_i$ has rank one and constant diagonal $1/v$, there is a sign
vector $\eps\in\{\pm1\}^X$ such that
$E_i=\eps\eps^{\mathsf T}/v$.
Put $D=\diag(\eps)$. Then $v(E_i\circ M)=DMD$ for every matrix $M$.
Schur closure shows that conjugation by $D$ preserves the
Bose--Mesner algebra. Since $D^2=I$, this map is an involutory
algebra automorphism and permutes its primitive idempotents.
Conjugation preserves ranks, and $DE_0D=E_i$.
\end{proof}

For singleton factors, put
$\mathcal K(i,j)=\{k:q_{ij}^k>0\}$.
~\Cref{thm:coefficient-criteria} says that
$E_i\circ E_j$ is projection-valued exactly when the positive
coefficients $q_{ij}^k$ are all equal. For a $P$-polynomial scheme with $d\ge2$ and $i,j\ne0$, the condition
$|\mathcal K(i,j)|=1$ is equivalent to equality in Pascasio's
inequality \cite[Corollary~4.2]{pascasio-inequality}.
Allowing more than one output idempotent requires equality of all
the coefficients on that support.

\Needspace{6\baselineskip}
In a $Q$-polynomial ordering, the only possible nonzero coefficients
in $E_1\circ E_i$ are
\[
 q_{1i}^{i-1},\qquad q_{1i}^{i},\qquad q_{1i}^{i+1},
\]
with indices outside $\{0,\ldots,d\}$ omitted. Thus the criterion
reduces to equality of the positive entries in this list.
For example, if $q_{1d}^d=0$, then $E_1\circ E_d$ is a scalar
multiple of $E_{d-1}$; the normalized identity
$v(E_1\circ E_d)=E_{d-1}$ holds exactly when
$q_{1d}^{d-1}=1$. These are direct specializations of the coefficient
criterion and do not require separate existence arguments.
Finally, factorizations are preserved under direct products, since
\begin{equation}\label{eq:product-construction}
 (M_1\otimes N_1)\circ(M_2\otimes N_2)
 =(M_1\circ M_2)\otimes(N_1\circ N_2).
\end{equation}
If the two input identities have parameters $\mu_1,\mu_2$, the
product-scheme identity has parameter $\mu_1\mu_2$.

\Cref{sec:dual-near} develops dual near-factorizations.
\Cref{sec:examples} gives examples with full and proper nonprincipal output supports.

\section{Dual near-factorizations}\label{sec:dual-near}

We consider the case in which the output contains every nonprincipal
idempotent.  This is the dual form of the near-factorization problem
for unions of relations studied in \cite{herman2026near}.

\begin{defn}\label{def:dual-near}
Let $S,T$ be nonempty subsets of $\{0,1,\ldots,d\}$.
We call $(S,T)$ a \defin{dual near-factorization with parameter $\mu$}
if $\mu>0$ and
\begin{equation}\label{eq:dual-near}
 v(E_S\circ E_T)=\mu(I-E_0).
\end{equation}
\end{defn}

The normalization in \eqref{eq:dual-near} makes $\mu$ an aggregate
Krein coefficient.  Suppose that a scheme with adjacency matrices
$\widehat A_0,\ldots,\widehat A_d$ has intersection numbers
$\widehat p_{ij}^{\,k}=q_{ij}^{\,k}$.  Then \eqref{eq:dual-near}
is equivalent to
\[
 \widehat A_S\widehat A_T=\mu(J-I),
 \qquad \widehat A_S=\sum_{i\in S}\widehat A_i.
\]
Such a combinatorial dual need not exist.  If it does, then $\mu$
is a positive integer.  The same integrality condition holds for
formally self-dual schemes.  The pair
$\bigl(\{0\},\{1,\ldots,d\}\bigr)$ is a trivial example with
$\mu=1$.  For the rest of this section we assume
$0\notin S\cup T$, so that both factors give spherical $2$-designs.

Put $r=m(S)$ and $s=m(T)$.  The principal coefficient, diagonal
entries, and rank inequality in \eqref{eq:dual-near} give
\begin{equation}\label{eq:dual-near-bounds}
 S\cap T=\varnothing,\qquad
 \mu=\frac{rs}{v-1},\qquad r+s\le v-1\le rs.
\end{equation}
Indeed, the principal coefficient is
$\sum_{i\in S\cap T}m_i$, and
$\rank(E_S\circ E_T)\le\rank(E_S)\rank(E_T)$.
In particular, $r,s\ge2$, $v\ge5$, and
\[
 1\le\mu\le\frac{v-1}{4}.
\]
The upper equality holds precisely when $r=s=(v-1)/2$, since
$rs\le(r+s)^2/4\le(v-1)^2/4$.
Moreover,
\begin{equation}\label{eq:dual-near-simplex}
 \left(\frac vrE_S\right)\circ
 \left(\frac vsE_T\right)=\frac{v}{v-1}(I-E_0).
\end{equation}
Thus every dual near-factorization gives a regular simplex after
normalization. Its strength is exactly two, since its cubic Gram
moment is $v[1-(v-1)^{-2}]>0$. The parameter $\mu=1$ is precisely the case in
which its rank is the product of the two input ranks.

\subsection{A criterion for a fixed factor}

Every entry of $I-E_0$ is nonzero.  Hence a necessary condition
for \eqref{eq:dual-near} is
$\eta_h(T)\ne0$ for $1\le h\le d$.
Under this condition, define
\begin{equation}\label{eq:reciprocal-matrix}
 C_T=\frac{v-1}{m(T)}I-
       \sum_{h=1}^d\frac{1}{\eta_h(T)}A_h.
\end{equation}

\begin{prop}\label{prop:reciprocal-criterion}
Let $T\subseteq\{1,\ldots,d\}$ be nonempty, and suppose
$\eta_h(T)\ne0$ for $1\le h\le d$.
There exists a nonempty $S\subseteq\{1,\ldots,d\}$ satisfying
\eqref{eq:dual-near} if and only if $C_T\one=0$ and all nonzero
eigenvalues of $C_T$ are equal to one positive number $c$.
In that case,
\[
 E_S=\frac1cC_T,\qquad \mu=\frac vc,
\]
and both $S$ and $\mu$ are uniquely determined by $T$.
\end{prop}

\begin{proof}
The entries in \eqref{eq:reciprocal-matrix} give
$C_T\circ E_T=I-E_0$.  Thus \eqref{eq:dual-near} is equivalent
to $E_S=(\mu/v)C_T$.  Since $C_T$ belongs to the Bose--Mesner
algebra, $C_T/c$ is a sum of primitive idempotents if and only if
its eigenvalues belong to $\{0,1\}$.  The condition
$C_T\one=0$ excludes $E_0$, and
$\tr C_T=v(v-1)/m(T)>0$ excludes the zero matrix.
This proves the equivalence.  If two solutions exist, their
idempotent sums are positive scalar multiples of one another.
Two nonzero idempotents can be scalar multiples only when that
scalar is one.  This proves uniqueness.
\end{proof}

The eigenvalue of $C_T$ on $E_j$ is
\begin{equation}\label{eq:reciprocal-eigenvalues}
 \gamma_j(T)=\frac{v-1}{m(T)}-
       \sum_{h=1}^d\frac{P_{jh}}{\eta_h(T)}.
\end{equation}
Consequently, for fixed $T$ the criterion requires only the
eigenmatrix entries, without a search over $S$.
The condition $\gamma_0(T)=0$ gives the useful necessary identity
\begin{equation}\label{eq:reciprocal-row-sum}
 \sum_{h=1}^d\frac{k_h}{\eta_h(T)}
       =\frac{v-1}{m(T)}.
\end{equation}

\subsection{Complementary input supports}

\begin{thm}
\label{thm:complementary-simplex}
Let $\mathfrak X$ be a symmetric association scheme on $v$ vertices, and
let $S,T$ be a partition of its nonprincipal primitive idempotent indices,
with $S,T\ne\varnothing$. Put
\[
 E=E_S,\qquad F=E_T=I-E_0-E,\qquad
 r=\operatorname{rank}E,\qquad s=\operatorname{rank}F=v-1-r.
\]
The following conditions are equivalent.
\begin{enumerate}
\renewcommand{\labelenumi}{(\roman{enumi})}
\item There is a real number $\mu>0$ such that $v(E\circ F)=\mu(I-E_0)$.
\item The matrices $E_0,E,F$ are the primitive idempotents of a
two-class fusion of $\mathfrak X$ whose strongly regular graph has
parameters $(v,k,a,c)$ satisfying
$2k=v+a-c$.
\item The matrices $E_0,E,F$ are the primitive idempotents of a
two-class fusion of $\mathfrak X$, and the two relation graphs of this
fusion have the same number of common neighbors of nonadjacent vertices.
\end{enumerate}
When these conditions hold, the fusion is primitive,
\[
 r,s\geq 2,\qquad \mu=\frac{rs}{v-1},
\]
and the fusion is formally self-dual if and only if its graph has
conference parameters
\[
 \left(v,\frac{v-1}{2},\frac{v-5}{4},\frac{v-1}{4}\right).
\]
\end{thm}

\begin{proof}
Suppose (i) holds. Comparison of diagonal entries gives $\mu=rs/(v-1)$. 
For distinct vertices $x,y$, put $\eta=vE_{xy}$.
Since $vF_{xy}=-1-\eta$, the off-diagonal equation gives $\eta^2+\eta=\mu$.
Thus the off-diagonal entries of $vE$ belong to $\{u,w\}$, where
\[
 u=\frac{-1+\sqrt{1+4\mu}}2>0,\qquad
 w=\frac{-1-\sqrt{1+4\mu}}2<0,
 \qquad u+w=-1.
\]
Both values occur. Otherwise $E\mathbf 1=0$ would make $E$ a scalar
multiple of $I-E_0$, and idempotence would force $E=I-E_0$, contrary
to $F\ne0$.

Next we define $A_{xy}=1$ when $x\ne y$ and $vE_{xy}=u$, and set all other
entries of $A$ equal to zero. 
Then we obtain $vE=(r-w)I+(u-w)A+wJ$.
The equation $E\mathbf 1=0$ shows that $A$ is regular. 
Moreover, we have
\[
 \operatorname{span}\{I,J,A\}
   =\operatorname{span}\{E_0,E,F\}
\]
is closed under ordinary multiplication. Hence $A$ is strongly regular.
Since $E$ belongs to the Bose--Mesner algebra of $\mathfrak X$,
$A$ and $J-I-A$ are unions of relations of $\mathfrak X$.
They therefore form a two-class fusion with primitive idempotents
$E_0,E,F$.

To identify its parameters and prove the converse, consider any such
two-class fusion. Write $\theta>\tau$ for the restricted eigenvalues
of its graph. Interchanging $E,F$ if necessary, the spectral formula is
\[
 E=\frac{1}{\theta-\tau}
       \left(A-\tau I-\frac{k-\tau}{v}J\right).
\]
The two off-diagonal entries of $vE$ are consequently
\[
 u=\frac{v-k+\tau}{\theta-\tau},\qquad
 w=\frac{\tau-k}{\theta-\tau}.
\]
Their sum is $-1$ if and only if
$2k=v+\theta+\tau=v+a-c$. This proves (i)$\Rightarrow$(ii).
Conversely, this parameter equation implies $u+w=-1$, so the
off-diagonal entries of $E\circ F$ are constant, say $\gamma$.
Since $EF=0$, each row gives
\[
 0=(EF)_{xx}=\frac{rs}{v^2}+(v-1)\gamma.
\]
Together with the diagonal entries, this proves (i), with
$\mu=rs/(v-1)$.

The complement of a strongly regular graph with parameters
$(v,k,a,c)$ has nonadjacent common-neighbor parameter $v-2k+a$.
Thus (ii) and (iii) are equivalent. Under (ii), the usual parameter
equation $k(k-a-1)=(v-k-1)c$ gives
\[
 c=\frac{k(v-k-1)}{v-1}>0.
\]
Both relation graphs therefore have positive nonadjacent
common-neighbor parameter and are connected. Hence the fusion is
primitive.

The Schur product of Gram matrices of ranks $r,s$ has rank at most
$rs$, by taking tensor products of their representing vectors.
Consequently $r+s=v-1\leq rs$, which excludes $r=1$ and $s=1$.

It remains to determine formal self-duality. Put $k'=v-k-1$, and let
$f$ be the multiplicity of $\theta$. The equations
$k+f\theta+(v-1-f)\tau=0$ and $2k=v+\theta+\tau$ yield
\[
 f-k=\frac{-(\theta+\tau+1)(k-\tau)}{\theta-\tau},\qquad
 f-k'=\frac{-(\theta+\tau+1)(k-\theta)}{\theta-\tau}.
\]

The two nonprincipal valencies are $k,k'=v-k-1$, whereas the corresponding multiplicities are $f,v-1-f$. 
Formal self-duality requires these
two pairs to agree up to permutation, and hence $f\in\{k,k'\}$.

Primitivity gives $k>\theta>\tau$, so $\theta+\tau=-1$ and
$k=(v-1)/2$. The formula for $c$ then gives
$c=(v-1)/4$ and $a=(v-5)/4$.
Conversely, for these conference parameters, compatible labelings give
\[
 P=Q=
 \begin{pmatrix}
  1&(v-1)/2&(v-1)/2\\
  1&(-1+\sqrt v)/2&(-1-\sqrt v)/2\\
  1&(-1-\sqrt v)/2&(-1+\sqrt v)/2
 \end{pmatrix}.
\]
Thus the fusion is formally self-dual.
\end{proof}

\begin{rmk}
The normalized Gram product in the theorem is
\[
 \left(\frac vr E\right)\circ\left(\frac vs F\right)
     =\frac{v}{v-1}(I-E_0),
\]
so the output consists of the $v$ vertices of a regular simplex.
If $\mu=1$, then $r+s=rs=v-1$, whence
$(r-1)(s-1)=1$ and $r=s=2$, $v=5$.
The fusion is therefore the scheme of the $5$-cycle.
This is the coefficient-one counterpart of
Herman and Miraftab's complementary-support theorem
\cite[Theorem~8.1]{herman2026association}, which assumes that $S,T$
partition the nonprincipal relation indices and
$A_SA_T=J-I$. No formal self-duality of the ambient scheme is
assumed in the theorem above.
\end{rmk}

\subsection{The last idempotent of a Hamming scheme}

The existence direction of the following result is the formally
self-dual form of the ternary Hamming construction in
\cite[Theorem~5.13]{herman2026near}.  The reciprocal
row-sum condition also determines when it can occur with
nonprincipal factors.

\begin{thm}\label{thm:hamming-endpoint}
Let $n\ge1$ and $q\ge2$, and use the usual idempotent ordering of $H(n,q)$.  
There exist a nonempty
$S\subseteq\{1,\ldots,n\}$ and $\mu>0$ such that $q^n(E_S\circ E_n)=\mu(I-E_0)$
if and only if $q=3$ and $n$ is even.  
In that case the unique
solution is
\begin{equation}\label{eq:ternary-dual-near}
 S=\{1,3,\ldots,n-1\},\qquad \mu=2^{n-1}.
\end{equation}
\end{thm}

\begin{proof}
The Hamming eigenmatrix formulas give
\[
 m_n=(q-1)^n,\qquad k_h=\binom nh(q-1)^h,
 \qquad Q_{hn}=(-1)^h(q-1)^{n-h}.
\]
Substituting these expressions in
\eqref{eq:reciprocal-row-sum} gives
\[
 \frac{\bigl(1-(q-1)^2\bigr)^n-1}{(q-1)^n}
       =\frac{q^n-1}{(q-1)^n}.
\]
Since $1-(q-1)^2=q(2-q)$, this is equivalent to $(2-q)^n=1$, which holds precisely when $q=3$ and $n$ is even.

For completeness, let $n$ be even and realize $H(n,3)$ on
$\F_3^n$.  
Let us write $F_a$ for the rank-one Fourier idempotent indexed by $a\in\F_3^n$.  
Then one can see that  $3^n(F_a\circ F_b)=F_{a+b}$.  
Next we put
\[
 A=\{a:\operatorname{wt}(a)\text{ is odd}\},\qquad
 B=(\F_3\setminus\{0\})^n.
\]
There is no representation $\boldsymbol0=a+b$ with
$a\in A$ and $b\in B$.  If $c\ne\boldsymbol0$ has weight $h$,
there are two choices of $b_i$ at each of its $n-h$ zero
coordinates.  At each nonzero coordinate, one choice gives
$a_i=0$ and the other gives $a_i\ne0$.  Exactly half of these
$2^h$ choices make $\operatorname{wt}(a)$ odd.
Thus every nonzero $c$ has $2^{n-1}$ representations.
Summing the Fourier idempotents proves
\eqref{eq:ternary-dual-near}.  
Uniqueness follows from~\Cref{prop:reciprocal-criterion}.
\end{proof}

For odd $n$, the same parity argument uses the even-weight idempotents, including $E_0$.  
This explains the role of the nonprincipal restriction in~\Cref{thm:hamming-endpoint}.

\section{Examples}\label{sec:examples}

We group the examples according to whether $E_U=I-E_0$ or
$E_U$ has proper nonprincipal support.

\subsection{Dual near-factorizations}\label{subsec:dual-near-examples}

The next example is about the Delsarte graph.
\begin{exa}\label{ex:delsarte-graph}
Let $\mathfrak X$ be the two-class scheme of the Delsarte graph,
with parameters $$(243,110,37,60),$$ 
see \cite{delsarte,BCN}.
Its second eigenmatrix is
\[
 Q=\begin{pmatrix}
 1&220&22\\
 1&4&-5\\
 1&-5&4
 \end{pmatrix}.
\]
The product of its nonprincipal columns is $20(242,-1,-1)$,
and therefore
\begin{equation}\label{eq:delsarte-product}
 243(E_1\circ E_2)=20(E_1+E_2).
\end{equation}
The valencies $110,132$ differ from the multiplicities $220,22$,
so the scheme is not formally self-dual.

The same scheme has no nontrivial identity
$A_SA_T=\rho(J-I)$, even when $A_0$ is allowed in a factor.
Such an identity forces $S\cap T=\varnothing$ by its diagonal
entries.  A factor $I$ gives only the trivial pair $I,J-I$.
The remaining disjoint pairs, up to interchange, give
\[
 \begin{aligned}
 A_1A_2&=72A_1+50A_2,\\
 (I+A_1)A_2&=72A_1+51A_2,\\
 A_1(I+A_2)&=73A_1+50A_2.
 \end{aligned}
\]
None is a scalar multiple of $J-I$.
Thus \eqref{eq:delsarte-product} is not an ordinary
near-factorization in the same scheme.
\end{exa}
The next example is about the Johnson scheme $J(9,4)$.
\begin{exa}\label{ex:johnson-nine-four}
Let $E_0,\ldots,E_4$ be the primitive idempotents of $J(9,4)$
in the usual ordering.  The valencies are $(1,20,60,40,5)$ and
the multiplicities are $(1,8,27,48,42)$; see \cite{BCN}.
The second eigenmatrix gives
\[
 \begin{split}
 (\eta_h(\{2,3\}))_{h=0}^4&=(75,3,-6,6,-3),\\
 (\eta_h(\{4\}))_{h=0}^4&=
 \left(42,-\frac{42}{5},\frac{21}{5},
                   -\frac{21}{5},\frac{42}{5}\right).
 \end{split}
\]
Their entrywise product is
$\frac{126}{5}(125,-1,-1,-1,-1)$.  Hence
\begin{equation}\label{eq:johnson-nine-four}
 126\bigl((E_2+E_3)\circ E_4\bigr)
       =\frac{126}{5}(I-E_0).
\end{equation}
This example has a sum of two primitive idempotents in one
factor and a nonintegral parameter.  Consequently, its Krein
parameters cannot be the intersection numbers of a combinatorial
dual.  In particular, the scheme is not formally self-dual.
\end{exa}

In the next example, we discuss the aggregation in $H(4,3)$.
\begin{exa}\label{ex:hamming-aggregation}
~\Cref{thm:hamming-endpoint} gives
$81\bigl((E_1+E_3)\circ E_4\bigr)=8(I-E_0)$, and $(r,s,t)=(40,16,80)$.
The constituent Krein vectors on output indices $1,2,3,4$ are
\[
 (q_{1,4}^k)_{k=1}^4=(0,0,2,4),\qquad
 (q_{3,4}^k)_{k=1}^4=(8,8,6,4).
\]
Their sum is constant, but neither is constant on its positive support.
Thus neither constituent product is projection-valued. The unused
idempotent $E_2$ has rank $24$, so full output does not require
complementary input supports.

More generally, for even $n\ge4$ in the ternary construction, none of
$E_i\circ E_n$ with $i$ odd is projection-valued. Indeed,
\[
 q_{i,n}^h=2^{n-h}\binom{h}{i-n+h}
 \qquad(n-i\le h\le n),
\]
and these are all its positive coefficients. Equality of the last two
would force $2i=n$. For $i=1$ this gives $n=2$, and for $i\ge2$ the ratio
of the third-last to the last coefficient would then be
\[
 \frac{4\binom{n-2}{i-2}}{\binom ni}=\frac{n-2}{n-1}\ne1.
\]
The normalized $E_n$ embedding is the tensor product of $n$ equilateral
triangles. The existence of the identity is transferred from
\cite[Theorem~5.13]{herman2026near}; the calculation above records why
aggregation is necessary in this realization.
\end{exa}

\begin{rmk}\label{rem:other-dual-near}
Several familiar examples fit the same definition.  In the
cycle scheme on $4a+1$ vertices, label $E_j$ by the Fourier
pair $\{j,-j\}$.  Then
\[
 (4a+1)(E_S\circ E_1)=I-E_0,
 \qquad
 S=\{j:1\le j\le2a,\ j\equiv2,3\pmod4\}.
\]
Indeed, for
$A=\{4j+2,4j+3:0\le j<a\}\subset\Z_{4a+1}$,
the translates $A-1$ and $A+1$ partition the nonzero residues.
Conference graphs give
$v(E_1\circ E_2)=\frac{v-1}{4}(I-E_0)$ by
~\Cref{thm:complementary-simplex}.
The order-four cyclotomic scheme on $\F_{17}$ also has
parameter one: for $H=\{\pm1,\pm4\}$, the sixteen sums in
$H+2H$ are precisely the nonzero elements of $\F_{17}$.
These schemes are formally self-dual, and the displayed identities
are dual forms of the known ordinary constructions; see
\cite[Theorem~3.1(iii), Lemma~9.2 and Proposition~10.4]{herman2026near}.
\end{rmk}

\subsection{Factorizations with proper nonprincipal output}
\label{subsec:proper-output-examples}

For $n\ge2k\ge6$, let $\overline X_{n,k}$ be the quotient of the
nonbinary Johnson scheme $J_3(n,k)$ obtained by identifying each
word of weight $k$ over $\{0,\pm1\}$ with its negative. Its vertices
are pairs $(S,\eps)$ with $S\in\tbinom{[n]}k$ and
$\eps\in\{\pm1\}^S$, modulo $(S,\eps)\sim(S,-\eps)$, and
$v=2^{k-1}\binom nk$.
The primitive idempotents $E_{a,b}$ of the quotient are indexed by
even $a$, $0\le a\le k$, and $0\le b\le k-a$, with
\[
 h_b(N)=\binom Nb-\binom N{b-1},\qquad
 m_{a,b}=\binom na h_b(n-a).
\]
Here $\binom N{-1}=0$. The eigenmatrix and multiplicity formulas used here follow from Tarnanen, Aaltonen and Goethals \cite{tarnanen1985nonbinary}. 
Passing to the quotient retains
precisely the even sign degrees.
In fact the sign reversal acts on the
sign-degree-$a$ component by multiplication by $(-1)^a$. 
The quotient retains exactly
the invariant components, namely those
with even $a$.

\begin{thm}\label{thm:signed-johnson}
Let $j$ be even with $2\le j<k$, and put
$E_L=E_{j,0}$, $E_B=E_{0,k}$ and $E_H=E_{j,k-j}$.
Then
\begin{equation}\label{eq:signed-johnson-product}
 v(E_L\circ E_B)=\mu E_H,\qquad
 \mu=\frac{h_k(n)}{h_{k-j}(n-j)},
\end{equation}
with factor and output ranks
\begin{equation}\label{eq:signed-johnson-ranks}
 (r,s,t)=\left(\binom nj,h_k(n),
                    \binom njh_{k-j}(n-j)\right).
\end{equation}
The normalized output consists of $v$ distinct nonantipodal points.
If $3j>2k$, it has strength exactly three.
\end{thm}

\begin{proof}
The normalized kernel of the last Johnson idempotent is
\[
 F_{N,a}(h)=\frac{(-1)^{a-h}}{\binom{N-a}{a-h}},\qquad
 F_{n-j,k-j}(h-j)=F_{n,k}(h).
\]
For $h\ge j$, we have
\[
 F_{n-j,k-j}(h-j)
 =\frac{(-1)^{k-h}}{\binom{n-k}{k-h}}
 =F_{n,k}(h).
\]

For two vertices $(S,\eps),(T,\delta)$, put $h=|S\cap T|$ and
$B_j=\binom kj^{-1}\sum_{J\in\binom{S\cap T}j}\eps_J\delta_J$,
where $\eps_J=\prod_{i\in J}\eps_i$.
The dual eigenmatrix formula gives normalized Gram entries.

Specializing the nonbinary Johnson eigenmatrix formulas of \cite{tarnanen1985nonbinary} to the even sign degrees of the quotient gives the following normalized Gram entries:
\[
 K_L=B_j,\qquad K_B=F_{n,k}(h),\qquad
 K_H=B_jF_{n-j,k-j}(h-j),
\]
where $K_H=0$ if $h<j$.

When $h<j$, we have $B_j=0$, and the kernel
of $E_{j,k-j}$ also vanishes. Therefore
$K_L\circ K_B=K_H$ for every pair of vertices,
proving the factorization after scaling.

Explicitly, the output vectors can be written as
\[
 z_{S,\eps}=\binom kj^{-1/2}
       \bigoplus_{J\in\binom Sj}\eps_J
       u_{S\setminus J}^{[n]\setminus J,k-j},
\]
where $u_A^{Y,a}$ is the unit evaluation vector of the last idempotent of $J(|Y|,a)$, indexed by $A\subseteq Y$; different $J$ use orthogonal blocks.
The nonzero blocks determine $S$.  Within a fixed support, equality or opposition requires all products $\prod_{i\in J}(\eps_i\delta_i)$ to have the same prescribed sign.
Comparing $j$-subsets differing in one element forces all $\eps_i\delta_i$ to agree.  
Since $j$ is even, only equality occurs, precisely for the same sign class.

The output is a spherical $2$-design because $E_H$ is nonprincipal.
A cubic coordinate monomial has zero sign average unless $J_1\mathbin{\triangle}J_2\mathbin{\triangle}J_3=\varnothing$.
For fixed $S$, averaging over sign classes is equivalent
to averaging over all sign vectors, since $j$ is even.
The sign factor of a cubic coordinate monomial is
$\eps_{J_1}\eps_{J_2}\eps_{J_3}
 =\eps_{J_1\mathbin{\triangle}J_2
              \mathbin{\triangle}J_3}$.
Its average vanishes unless
$J_1\mathbin{\triangle}J_2
\mathbin{\triangle}J_3=\varnothing$.
In that case, every element of
$J_1\cup J_2\cup J_3$ belongs to exactly
two of the three sets. Consequently,
\[
 |J_1\cup J_2\cup J_3|=\frac{3j}{2}\le k.
\]

Hence $3j>2k$ gives strength three.  Both input ranks exceed one,
so~\Cref{prop:no-tensor-4design} excludes strength four.
\end{proof}

\begin{prop}\label{prop:cyclotomic-product}
Let $p>81$ be a prime with $p\equiv1\pmod4$.  
We choose
$i\in\F_p$ with $i^2=-1$, and put
$H=\{1,-1,i,-i\}\le\F_p^\times$.
In the symmetric cyclotomic scheme of $\F_p$, denote by $E_{cH}$
the primitive idempotent indexed by the multiplicative coset $cH$.
Then
\begin{equation}\label{eq:cyclotomic-product}
 p(E_H\circ E_{2H})
 =E_H+E_{3H}+E_{(1+2i)H}+E_{(1-2i)H}.
\end{equation}
The normalized Gram-matrix product is a nonantipodal spherical $3$-design with $p$ points in $\R^{16}$.
In particular, there are infinitely many primitive association schemes admitting this construction.
\end{prop}

\begin{proof}
For the additive characters of $\F_p$, let $F_a$ be the corresponding rank-one Fourier idempotents.  
They satisfy $p(F_a\circ F_b)=F_{a+b}$.  
Thus the Krein coefficients in the
cyclotomic scheme are the multiplicities of sums from two cosets.
In $\Z[\mathrm i]$, the sixteen elements of
\[
 \mathcal U=\{\pm1,\pm\mathrm i\}
                 +2\{\pm1,\pm\mathrm i\}
\]
are distinct and form the four orbits represented by $1,3,1+2\mathrm i,1-2\mathrm i$ under multiplication by $\{\pm1,\pm\mathrm i\}$.  
To justify reduction modulo $p$, observe that if a Gaussian integer $z=a+b\mathrm i$ maps to zero in $\F_p$, then we have 
$a^2+b^2=(a+bi)(a-bi)\equiv0\pmod p$.
Therefore a nonzero Gaussian integer with
squared absolute value less than $p$ cannot reduce to zero.

The difference of two distinct elements of $\mathcal U$ has squared absolute value at most $36$, so its
reduction in $\F_p$ remains nonzero.  This proves \eqref{eq:cyclotomic-product}.

Every element of $\mathcal U$ has odd coordinate sum.  
Hence a sum of three such elements cannot vanish in $\Z[\mathrm i]$.
Such a sum has squared absolute value at most $81$, so it also remains nonzero modulo $p$.  
Thus its reduction $U\subset\F_p$
satisfies $(U+U)\cap U=\varnothing$.  
Since $U=-U$, this means that no three
elements of $U$ sum to zero. Equivalently,
the aggregate Krein coefficients of the
output vanish on its own support, which
is precisely the condition in
\eqref{eq:scheme-3design-criterion}.

The output embedding is therefore a spherical $3$-design by~\Cref{prop:scheme-3design}.
The output points are distinct, since every nontrivial additive character of $\F_p$ is faithful.
It has no antipodal pair: an additive character of the odd-order
group $\F_p$ cannot take the value $-1$.
The scheme is primitive since the additive group of $\F_p$ has
prime order.  
Since Dirichlet's theorem on primes 
guarantees infinitely many primes $p\equiv 1\pmod 4$,
there are infinitely many primitive association schemes admitting this construction.
\end{proof}

The next example is about the binary Hamming schemes.
\begin{exa}\label{ex:hamming-product}
Let $E_0,\ldots,E_{2m}$ be the primitive idempotents of the binary
Hamming scheme $H(2m,2)$, where $m\ge2$.  Counting pairs of
binary characters by their weights yields
\begin{equation}\label{eq:hamming-product}
 2^{2m}(E_1\circ E_m)
 =(m+1)(E_{m-1}+E_{m+1}).
\end{equation}
In particular, we have 
$16(E_1\circ E_2)=3(E_1+E_3)$ in $H(4,2)$.
The corresponding tensor-product spherical $2$-design is the
$8$-dimensional cross-polytope, hence a spherical $3$-design.
More generally, the output in \eqref{eq:hamming-product} is a
spherical $3$-design whenever $m$ is even, since the indices
$m-1$ and $m+1$ are odd.
In contrast, $8(E_1\circ E_2)=2E_1+3E_3$ in $H(3,2)$,
so orthogonality of the input idempotents alone does not suffice.
\end{exa}

\section*{Data Availability} No datasets were generated or analysed during the current study.
\section*{Declarations}
The authors declare no conflict of interest.

\bibliographystyle{plainnat}
\bibliography{MPF}
\newpage
    
\end{document}